\documentclass[11pt,leqno]{article}

\usepackage{vmargin}

\usepackage[latin1]{inputenc}
\usepackage[english]{babel}
\usepackage{amsmath}
\usepackage{amssymb}

\usepackage{tabularx}
\usepackage{epsfig, floatflt,amssymb}
\usepackage{moreverb} %% pour le verbatim en boite
\usepackage{multirow} %% pour regrouper un texte sur plusieurs lignes dans une table
\usepackage{url} %% pour citer les url par \url
\usepackage[all]{xy} %% pour la barre au dessus des symboles
\usepackage{shorttoc} %% pour plusieurs tables des matières par la commande \shorttableofcontents{Titre}{profondeur}.
\usepackage{textcomp} %% pour le symbol pour mille par \textperthousand.
\usepackage{float}
\usepackage[right]{eurosym}

\newcommand{\diver}{\operatorname{div}}

\makeindex %% crée l'index

\usepackage{latexsym}
\usepackage{amssymb,amsbsy,amsmath,amsfonts,amscd}
\usepackage{vmargin}
\usepackage{epsfig}
\usepackage{amsmath}
\usepackage{amsfonts}
\usepackage{amssymb}
\usepackage{graphicx}%
\usepackage[dvips]{color}

\newtheorem{lemme}{Lemma}[section]%a enlever peut etre danla la these voir these.sty

\newtheorem{Theorem}{Theorem}[section]
\newtheorem{Definition}{Definition}[section]

\newtheorem{Remark}{Remark}[section]

\newcommand{\R}{\ifmmode{{\rm I} \hskip -2pt {\rm R}}
    \else{\hbox{$I\hskip -2pt R$}}\fi}

\newcommand{\RR}{\mathbb{R}}

\newcommand{\tore}{\mathbb{T}_3}

\newcommand{\BEQ} {\begin{equation} }
\newcommand{\EEQ} {\end{equation} }
\makeatletter
\@addtoreset{equation}{section}

\usepackage{epsfig}

\newcounter{taskcounter}[section]

\newcounter{technique}[section]

\def\vec#1{\boldsymbol{#1}}

\def\diver{\mathop{\mathrm{div}}\nolimits}

\def\R{\mathcal{R}}

\def\bef{\vec{f}}

\def\bv{\vec{v}}

\def\b0{\vec{0}}

\def\bw{\vec{w}}

\newcommand{\diameter}{\operatorname{diameter}}

\begin{document}%

\begin{titlepage}
\title{\vskip-2cm\textbf{Singularity \& Regularity Issues for  Simplified  Models of Turbulence}\vskip0.5cm}

\author{\begin{Large}$\hbox{Hani Ali$^\ast$, Zied Ammari}\thanks{IRMAR, UMR 6625,
Universit\'e Rennes 1,
Campus Beaulieu,
35042 Rennes cedex
FRANCE. Email: hani.ali@univ-rennes1.fr, zied.ammari@univ-rennes1.fr}$ \end{Large}}
\end{titlepage}
\date{}
\maketitle

\begin{center}
 \textbf{Abstract}
\end{center}
\ \ We consider a  family of  Leray-$\alpha$ models with periodic boundary conditions in three space dimensions. Such models are a regularization, with respect to a parameter $\theta$,
of the Navier-Stokes equations. In particular, they share with the original equation (NS) the property of existence of global weak solutions.
We establish an upper bound on the Hausdorff dimension of the time singular set of those weak solutions when $\theta$  is subcritical. %(i.e., $0<\theta<1/4$).
The result is an  interpolation between the bound proved by Scheffer for the Navier-Stokes equations and the regularity result proved in  \cite{A01}.
\ \\

{\small\it MS Classification:} 35Q30, 35Q35, 76F60

\medskip


\section{Introduction}
We consider, for $\alpha>0$ and $0<\theta<1/4$,  the Leray-$\alpha$ equations in the 3-dimensional flat Torus $\mathbb{T}_3$
\begin{equation}
\label{alpha ns}
 \left\{
\begin{array} {llll} \displaystyle
 \frac{\partial \vec{u}^{}}{\partial t}+(\overline{\vec{u}^{}} \cdot \nabla)
\vec{u}^{} - \nu \Delta \vec{u}^{} + \nabla
p^{} = \vec{f} \ \ \ \ \ \hbox{ in }\ \RR^{+}\times\mathbb{T}_3,\\
\overline{\vec{u}}=(1-\alpha^{2} \Delta)^{{-\theta}{}}\vec{u},\ \ \hbox{ in } \ \mathbb{T}_3,  \\
\nabla \cdot \vec{u}^{}=0,\
 \displaystyle \int_{\mathbb{T}_3} \vec{u}^{}=0,\\
\vec{u}^{}_{|t=0}=\vec{u}_{0}^{}.
\end{array}\right.
\end{equation}
Here the unknowns are  the  velocity vector field $\vec{u}$ and the scalar pressure  $p$. The viscosity $\nu$, the initial velocity vector field $\vec{u}_{0}$ and the external
force $\vec{f}$, with $\nabla \cdot \vec{f} = 0$, are given. \\
The nonlocal operator $M_\theta=(1-\alpha^2 \Delta)^{{-\theta}{}}$, acting on $L^2(\mathbb{T}_3,\mathbb{R}^3)$,  is defined through the Fourier transform on the torus
\begin{equation}
\widehat{M_\theta \vec{u}}({\vec{k}})=(1+\alpha^2 |{\vec{k}}|^{2})^{-\theta} \, \widehat{\vec{u}}({\vec{k}})\;, \quad \vec{k}\in\mathbb{Z}^3.
\end{equation}

\ \ The Leray-$\alpha$ equations 
are among the simplest models of turbulence, introduced nearly a decade ago for numerical simulation purposes. 
When $\theta = 0$, (\ref{alpha ns}) reduces to the Navier Stokes equation for an incompressible fluid. For $\theta>0$, $\overline{\vec{u}}$ is a regularization of the 
velocity vector field $\vec{u}$. Actually, a crude regularization (or filtering) appeared in the early work of J.~Leray \cite{JL34} where a 
mollifier was used (i.e., $\overline{\vec{u}}=\phi_\varepsilon*\vec{u}$) instead of the operator $M_\theta$.\\

\ \  The Leray-$\alpha$ models  are  
approximation of the Navier-Stokes equations and in fact they have  several properties in common with (NSE). In particular, \eqref{alpha ns} have existence of weak solutions for arbitrary 
time and large initial-data (see Theorem \ref{1deuxieme}).     \\

\ \ Our goal, in this short note, is to establish an upper bound for the Hausdorff dimension of the time singular set $\mathcal{S}_\theta(\vec{u})$ of weak solutions $\vec{u}$ of \eqref{alpha ns}. We know, thanks  to 
Scheffer's work \cite{S76,S77-2}, that  if $\vec{u}$ is a weak Leray solution of the Navier-Stokes equations then the 
$\frac{1}{2}$-dimensional Hausdorff measure of the time singular set of $\vec{u}$  is zero. Further, 
when $\theta=\frac{1}{4}$, the author in \cite{A01} proved the existence  of a unique  regular weak solution to the Leray-$\alpha$ model \eqref{alpha ns}.
Therefore, it is intersecting to understand how the potential time singular set $\mathcal{S}_\theta(\vec{u})$  may depend on the regularization 
parameter $\theta$. \\

\ \ In fact, we will prove that $\frac{1-4\theta}{2}$-dimensional Hausdorff measure of  the time singular set $\mathcal{S}_\theta(\vec{u})$ of any 
weak solution $\vec{u}$ of  \eqref{alpha ns} is zero (see Theorem \ref{main}). \\

\ \  Although we  consider only the  Leray-$\alpha$ equations \eqref{alpha ns}, the same results hold for other models of turbulence as the 
magnetohydrodynamics MHD-$\alpha$ equations. It was  observed in \cite{HLT10} that the qualitative properties of the equation \eqref{alpha ns} 
relies on the regularization effect of the operator $M_\theta$ rather than its explicit form. Indeed, this fact can be checked for the Hausdorff dimension of the 
time singular set $\mathcal{S}_\theta(\vec{u})$.

\section{Preliminaries}
Before giving some preliminary results  we fix some notations. For $p \in [1 , \infty)$, the Lebesgue spaces $L^p(\tore)$, the Sobolev spaces $W^{1,p}(\tore)$
and the Bochner spaces $L^p (0, T ; X)$, $C(0, T ;$ $ X)$,  $X$ being a Banach space, are defined in a standard way. In
addition for $s \ge -1$, we introduce the spaces
$$\vec{V}^s=\left\{   \vec{u} \in  W^{s,2}(\tore)^3,\ \int_{\tore} \vec{u}=0,\ \diver \vec{u}=0 \right\},$$
endowed with the norms  
$$ \|\vec{u}\|_{\vec{V}^s}^2=\sum_{\vec{k}\in\mathbb{Z}^3} |\vec{k}|^{2s} \,|\widehat{\vec{u}}({\vec{k}})|^2\,.$$
For the sake of simplicity we introduce the notations 
$$\vec{H}:=\vec{V}^0 \quad \hbox{ and } \quad \vec{V}:=\vec{V}^1\,.$$

\subsection{A priori estimates }
The essential feature of the operator $M_\theta$ is the following regularization effect.
\begin{lemme}
\label{regularitedeuxieme}
Let $\theta\in \RR_+,s \ge -1 $ and assume that $\vec{u} \in \vec{V}^{{s}}$.   Then $ M_\theta\vec{u} \in
\vec{V}^{{s}+ 2\theta}$  and
$$\|M_\theta\vec{u}  \|_{\vec{V}^{{s}+2\theta} }  \le \frac{1}{\alpha^{2\theta}} \|\vec{u} \|_{\vec{V}^{s}}.
  $$
\end{lemme}
\ \medskip
\ \ Next, we prove a priori estimates in the same manner as for the Navier Stokes equations (see \cite{RT83}). 
We suppose that $\vec{u}$ is a sufficient regular solution of \eqref{alpha ns}. 

\subsubsection{A priori estimates in $\vec{H}$}
\begin{lemme}
\label{l2deuxieme}
Let $\vec{f} \in L^{2}([0,T],\vec{V}^{{-1}})$ and $\vec{u}_{0} \in \vec{H} $, for all $T \ge 0$ there exists  $K_{1}(T)$ and $K_{2}(T)$  
such that any solution $\vec{u}$ of (\ref{alpha ns}) satisfies 
\begin{equation}
\|\vec{u}\|_{L^{2}([0,T],\vec{V})}^2 \le K_{1}(T), \ \ \  \hbox{where} \ K_{1}(T)= \frac{1}{\nu}\left( \|\vec{u}_{0}\|_{\vec{H}}^{2}+\frac{1}{\nu} \|\vec{f} \|_{L^{2}([0,T],\vec{V}^{{-1}})}^{2} \right)
\end{equation}
and
\begin{equation}
\|\vec{u}\|_{L^{\infty}([0,T],\vec{H})}^2 \le K_{2}(T) , \ \ \  \hbox{where} \  K_{2}(T)= {\nu} \, K_{1}(T)  .
\end{equation}
\end{lemme}
\textbf{Proof}
Taking the $L^2$-inner product of the first equation of (\ref{alpha ns})  with $\vec{u}$ and integrating by parts. Using the incompressibility of the velocity field and the duality relation we obtain
\begin{eqnarray*}
\frac12\frac{d}{dt}\|\vec{u}\|_{\vec{H}}^2+\nu \|\nabla \vec{u}\|_{\vec{H}}^2&=&\int_{\tore}\vec{f} \ \vec{u} dx
\le \|\vec{f}\|_{\vec{V}^{-1}} \|\vec{u}\|_{\vec{V}}.
\end{eqnarray*}
Using Young inequality we get
\begin{eqnarray*}
\frac{d}{dt}\|\vec{u}\|_{\vec{H}}^2+\nu\|\nabla \vec{u}\|_{\vec{H}}^2&\le&\frac{1}{\nu}\|\vec{f}\|_{{\vec{V}^{-1}}}^2.
\end{eqnarray*}
Integration with respect to  time gives the desired estimates. $\hfill\blacksquare$

\subsubsection{A priori estimates in $\vec{V}$ }
Now we use the regularization effect of  Lemma \ref{regularitedeuxieme} to prove the following a priori estimates.
\begin{lemme}
\label{l321}
Let $\vec{f} \in L^{2}([0,T], \vec{H})$ and $\vec{u}_{0} \in \vec{V}$ . Assume that $0<\theta< 1/4 $.
 Then there exists $T_*:= T_*(\vec{u}_{0})$ and  $M_{}(T_*)<\infty$  such that the solution $\vec{u}_{}^{}$ of \eqref{alpha ns}  satisfies
 $$ \sup_{t \in [0, T_*]} \|\vec{u}\|_{\vec{V}}^2 \le 2(1+ \|\vec{u}_{0} \|_{\vec{V}}^2)$$
 and
 $$ \int^{T_*}_{0}\|\Delta \vec{u}(t)\|_{\vec{H}}^2 dt \le M_{}(T_*).$$
 %where
\end{lemme}
\textbf{Proof}
Taking the $L^2$-inner product of the first equation of (\ref{alpha ns})  with $-\Delta \vec{u}$ and integrating by parts. Using  the incompressibility of the velocity field and the duality relation  
combined with H\"{o}lder inequality and Sobolev injection, we obtain
\begin{eqnarray*}
\frac12\frac{d}{dt}\|\nabla\vec{u}\|_{\vec{H}}^2+\nu\|\Delta \vec{u}\|_{\vec{H}}^2& \le & \int_{\tore}  |(\overline{\vec{u}}\cdot\nabla ) \vec{u} \Delta\vec{u}| dx + \int_{\tore}|\vec{f} \Delta\vec{u}|dx\\
&\le& \alpha^{-2\theta}\|\nabla{\vec{u}}\|_{\vec{H}}  \| \nabla \vec{u}\|_{\vec{V}^{\frac{1}{2}-2\theta}}   \|\Delta\vec{u}\|_{\vec{H}}+   \|\vec{f}\|_{\vec{H}} \| \Delta\vec{u}\|_{\vec{H}}.
\end{eqnarray*}
Interpolating between $\vec{V}^{1}$ and $\vec{V}^2$ we get
\begin{eqnarray*}
\frac12\frac{d}{dt}\|\nabla\vec{u}\|_{\vec{H}}^2+\nu \|\Delta \vec{u}\|_{\vec{H}}^2
&\le& \displaystyle \alpha^{-2\theta}\|\nabla{\vec{u}}\|_{\vec{H}}^{\frac{3}{2}+ 2 \theta}   \displaystyle \|\Delta\vec{u}\|_{\vec{H}}^{\frac{3}{2}- 2 \theta} +   \|\vec{f}\|_{\vec{H}} \|\Delta\vec{u}\|_{\vec{H}}.
\end{eqnarray*}
Using Young inequality we get
\begin{eqnarray}
\label{3.24}
\displaystyle \frac{d}{dt}\|\nabla \vec{u}\|_{\vec{H}}^2+ \displaystyle \nu\|\Delta \vec{u}\|_{\vec{H}}^2&\le& \displaystyle \frac{1}{\nu} \|\vec{f}\|_{\vec{H}}^2 + 
\displaystyle C(\alpha,\theta) \displaystyle\|\nabla \vec{u}\|_{\vec{H}}^{ \frac{2(3+4\theta)}{1+4 \theta}}.
\end{eqnarray}
 We get a differential inequality
 \begin{eqnarray}
\label{differentiel}
 \displaystyle Y^{'} \le \displaystyle C(\alpha,\theta, \nu , f ) \;Y^{\gamma},
 \end{eqnarray}
 where $$ Y(t)= 1+ \|\vec{u}\|_{\vec{V}}^2 \quad \hbox{ and }   \quad \gamma = \frac{3+4\theta}{1+4 \theta}\,.$$
 We conclude that
 $$ Y(t) \le \frac{Y(0)}{ ( 1-2Y(0)^{\gamma -1}   C(\alpha,\theta, \nu , f ) \,t  )^{\frac{1}{\gamma-1}}  }$$
 as long as $ \displaystyle t < \frac{1}{2Y(0)^{\gamma -1}  C(\alpha,\theta, \nu , f )}\,$. Hence we obtain
  $$ \sup_{t \in [0, T_*]} \|\vec{u}\|_{\vec{V}}^2 \le 2(1+ \|\vec{u}_{0}\|_{\vec{H}}^2)$$
  \begin{eqnarray}
  \label{maximal time}
  \hbox{ with  } \, \displaystyle T_{*}:= \frac{3}{8   C(\alpha,\theta, \nu , f )}\frac{1}{(1+\|\vec{u}_{0}\|_{\vec{V}}^2)^{\gamma -1} }.
  \end{eqnarray}
Integrating  (\ref{3.24}) with respect to  time on $[0, T_{*}] $  gives the desired estimates
 $$ \int^{T_*}_{0}\|\Delta \vec{u}(t)\|_{\vec{H}}^2 dt \le M_{}(T_*),$$
 where  $$ M_{}(T_*)= \frac{1}{\nu}\left(  \|\vec{u}_{0}\|_{\vec{H}}^2 + \frac{2}{\nu} \int^{T_*}_{0}\|\vec{f}\|_{\vec{H}}^2 dt +C(\alpha,\theta)\, [ 2(1+ \|\vec{u}_{0}\|_{\vec{H}}^2)]^{\gamma}\right). 
 $$ $\hfill\blacksquare$
 
\subsection{Existence and uniqueness results}
The next two theorems collect the most typical results for the Leray-$\alpha$ models of turbulence (see \cite{A01}, \cite{OT2007}).
The proofs of these two theorems follow by combination of the  above a priori estimates with a Galerkin method. This  is a classical argument which we avoid its repetition.   
For further information, we refer the reader to \cite{RT83}, \cite{A01} and the references therein. \\
 
By $\displaystyle C_{weak}([0,T];\vec{H})$ we denote the vector space of all mappings 
$\displaystyle \vec{v} : [0,T]\rightarrow \vec{H}$ such that for any $\displaystyle \vec{h} \in \vec{H},$ the function
$$ t \mapsto  \int_{\tore} \vec{v}(t) \vec{h} d\vec{x} $$
is continuous on $[0,T].$

\begin{Theorem}
\label{1deuxieme}
Let $\vec{f} \in L^{2}([0,T],\vec{V}^{-{{1}{}}})$ and $\vec{u}_{0} \in \vec{H}$. Assume that  $0\leq\theta< 1/4 $.
Then for any $T>0$  there exist a weak Leray solution $( \vec{u},p^{}):=( \vec{u}_{\alpha},p^{}_{\alpha})$ to (\ref{alpha ns})  such that
 $$  \vec{u}_{} \in  C_{weak}([0,T];\vec{H}) \cap L^{2}([0,T];\vec{V}), \quad \ \quad
 \displaystyle \frac{ \partial\vec{u}_{}^{}}{ \partial t} \in {L^{\frac{5}{3-2\theta}}([0,T];W^{-1,\frac{5}{3-2\theta} }(\mathbb{T}_3)^3)},$$
   $${p}_{} \in  L^{\frac{5}{3-2\theta}}([0,T],L^{{\frac{5}{3-2\theta}}}(\mathbb{T}_3)),$$
%Moreover, this solution depends continuously on the initial data $\vec{u}_{0}$ in the $L^{{2}}$ norm.
\begin{equation}
\begin{split}
\int_0^T \langle \frac{ \partial\vec{u}_{}^{}}{ \partial t} , \bw \rangle  -  (\overline{\vec{u}} \otimes \vec{u}, \nabla
\bw) +  \nu(
\nabla \bv, \nabla \bw )- (p, \diver \bw) \; dt\\
\qquad
= \int_0^T \langle \bef, \bw \rangle \; dt
\qquad \textrm{ for all } \bw\in L^{\frac{5}{2+2\theta}}(0,T; W^{1,\frac{5}{2+2\theta}}(\mathbb{T}_3)^3),
\end{split}\label{weak1001}
\end{equation}
where the velocity $\vec{u}$ verifies
\begin{equation}\begin{split}
\sup_{t\in(0,T)}\|{\vec{u}}(t)\|_{\vec{H}}^2+ \int_0^T \|{\vec{u}}(t)\|_{\vec{V}}^2 dt\le \|{\vec{u}}_0\|_{\vec{H}}^2 + \int_0^T \langle \bef,\vec{u} \rangle \; dt
%C(\bv_0,\bef) < \infty
, \label{regularity}
\end{split}
\end{equation}
and the initial data is attained in the following sense
   $$\lim_{t\to 0+} \|\vec{u}(t)-\vec{u}_0\|_{\vec{H}}^2 =0.$$
\end{Theorem}
\begin{Remark}
 If $ \theta = \frac{1}{4}$ a weak solution to Leray-$\alpha$ model is called regular weak solution \cite{A01}, in addition, the solution is unique and it satisfies
 $$  \vec{u}_{} \in  C([0,T];\vec{H}) \cap L^{2}([0,T];\vec{V}), \quad   \displaystyle \frac{ \partial\vec{u}_{}^{}}{ \partial t} \in L^{2}([0,T];\vec{V}^{-1}), \quad \hbox{ and } \quad
 {p}_{} \in  L^{2}([0,T],L^{2}(\mathbb{T}_3)).$$
 In this case one also has energy equality in (\ref{regularity}) instead of inequality.
\end{Remark}
\begin{Theorem}
\label{2deuxieme}
Let $\vec{f} \in L^{2}([0,T],\vec{H})$ and $\vec{u}_{0} \in \vec{V}$. Assume that $0\leq \theta< 1/4 $.
 Then there exists $T_*:= T_*(\vec{u}_{0})$, determined by (\ref{maximal time}), and  there exists a unique strong solution  $ \vec{u}$ to  
 (\ref{alpha ns}) on $[0, T_*]$ satisfying:
$$ \vec{u}_{} \in  C([0,T_*];\vec{V}) \cap L^{2}([0,T_*];\vec{V}^{2}),$$  $$ \displaystyle \frac{\partial\vec{u}_{}^{}}{\partial t} \in {L^{2}([0,T_*];L^{2}(\mathbb{T}_3)^3)} \quad \hbox{ and } \quad   
{p}_{} \in  L^{2}([0,T_*],W^{1,2}(\mathbb{T}_3))\,.$$\\
\end{Theorem}
\begin{Remark} 
If $ \theta = \frac{1}{4}$ the strong  solution to the Leray-$\alpha$ model exists for any arbitrary time $T >0$. 
Indeed,  when $ \theta =\frac{1}{4}$, $ \gamma = \frac{3+4\theta}{1+4 \theta} = 2$, the differential inequality \eqref{differentiel} becomes 
 \begin{eqnarray}
\label{differentiel 2}
 \displaystyle Y^{'} \le \displaystyle C(\alpha,\theta, \nu , f )\, Y^{2}\,.
 \end{eqnarray}
Thus, using Gronwall's inequality $(Y \in L^1 [0,T])$ we obtain the desired result.
\end{Remark}
\begin{Remark}
Let us assume that $\vec{f} \in L^{2}([0,T],\vec{H})$, $\vec{u}_{0} \in \vec{V}$  and  $0<\theta< 1/4 $.
With the local existence of strong solution and the weak=strong theorem of Serrin \cite{SE63} , the
solution $\vec{u}$ is automatically strong and thus smooth on $ [0,T_*)\times \tore $, where $T_* \in [0,T]$.  We
note that in  Theorem \ref{1deuxieme}, a weak solution satisfying the above properties with $T=\infty$ exists for
every divergence-free $\vec{u}_{0} \in \vec{H}$.
\end{Remark}

\section{The Main  Result  And Its Proof}
The basic facts about Hausdorff measures can be found for instance in \cite{Federer}. We recall here the definition of those measures.
\begin{Definition}
\label{definitionhaussdrof}
Let $X$ be a metric space and let $a > O$. The $a$-dimensional Hausdorff measure of a subset $Y$ of $X$ is
$$
\displaystyle
\mu_a(Y)=\lim_{\epsilon \searrow o}\mu_{a,\epsilon}(Y)=\sup_{\epsilon > 0}\mu_{a,\epsilon}(Y)
$$
where
$$
\displaystyle
\mu_{a,\epsilon}(Y)=\inf \sum_{j}(\diameter B_j)^a,
$$
the infimum  being taken over all the coverings of $Y$ by balls $B_j$ such that $\diameter B_j \le \epsilon$.
\end{Definition}
\begin{Definition}
Let $T>0$. We define  the time singular set $\mathcal{S}_\theta(\vec{u})$ of $\vec{u}(t) $, a weak solution of (\ref{alpha ns}) given by Theorem \ref{1deuxieme}, 
as the set of $t \in [0, T]$ such that $\vec{u}(t) \notin \vec{V}$.
\end{Definition}
The main result of the paper is the following theorem.
\begin{Theorem}
\label{main} Let $\vec{u}$ be any weak Leray solution to eqs. (\ref{alpha ns}) given by  Theorem \ref{1deuxieme} ( We suppose  that the external force $ \vec{f} \in L^{2}([0,T],\vec{H})$ ).
Then for any $T >0$ the $  \frac{1-4\theta}{2} $-dimensional Hausdorff measure of the time singular  set $\mathcal{S}_\theta(\vec{u})$ of $\vec{u}$   is zero.
\end{Theorem}
The rest of the paper is devoted to the proof of the main Theorem. 
The following Lemma  characterizes the structure of the time singularity set of a weak solution of (\ref{alpha ns}).
\begin{lemme}
\label{temam book}
We assume that $\vec{u}_{0} \in \vec{H}$,  $\vec{f} \in L^{2}([0,T],\vec{H})$ and $ \vec{u}$ is any weak solution of (\ref{alpha ns}) given by  Theorem \ref{1deuxieme}. 
Then there exist an open set $ \mathcal{O} $ of $(0,T) $  such that:\\
 (i) For all $t \in   \mathcal{O} $ there exist  $t \in (t_1, t_2)\subseteq (0,T) $ such that
$ \vec{u} \in C((t_1, t_2), \vec{V}) $.\\
(ii) The  Lebesgue measure of $[0,T]\setminus \mathcal{O} $  is zero.
%$(t_1, t_2) $ such that
\end{lemme}
\textbf{Proof.} Since $\vec{u}_{} \in  C_{weak}([0,T];\vec{H}) $, $\vec{u}_{}(t)$ is well defined for every $t$ and we can define
$$
\Sigma=\{ t \in [0,T], \vec{u}_{}(t) \in \vec{V} \},
$$
$$
\Sigma^c=\{ t \in [0,T], \vec{u}_{}(t) \notin \vec{V} \},
$$
$$
\mathcal{O}=\{ t \in (0,T),  \exists \epsilon > 0, \vec{u}_{}  \in  C((t-\epsilon, t+\epsilon),  \vec{V})  \}.
$$
It is clear that $\mathcal{O}$ is open.
Since $  \vec{u}_{}  \in  L^{2}([0,T];\vec{V})$, $\Sigma^c$ has Lebesgue measure zero. Let us take $t_0$ such  that  $ t_0  \in \Sigma, $  and $ t_0  \notin \mathcal{O} $, 
then according to Theorem \ref{2deuxieme}, there exists $ \epsilon >0$ such that $\vec{u}_{}  \in   C((t_0, t_0+\epsilon), \vec{V}) $. So that, $ t_0$  is 
the left end of one of the connected components of $\mathcal{O} $. Thus $\Sigma\setminus \mathcal{O} $  is countable and  consequently 
$ [0,T]\setminus \mathcal{O}$  has Lebesgue measure zero. This finishes the proof. $\hfill\blacksquare$\\
\begin{Remark}
We deduce from Theorem \ref{2deuxieme} that, if $(\alpha_i, \beta_i )$ , $i \in I$, is one of  the connected components of $ \mathcal{O}$, then
$$ \lim_{t\rightarrow \beta_i} \| \vec{u}(t)\|_{\vec{V}}=+\infty.$$
Indeed, otherwise Theorem \ref{2deuxieme} would show that there exist an $\epsilon >0$ such that $\vec{u}_{}  \in  C((\beta_i, \beta_i+\epsilon),  \vec{V})$ and $\beta_i$ 
would not be the end point of a connected  component of $ \mathcal{O}$.
\end{Remark}
\begin{lemme}
\label{temam book2}
Under the same notations of Lemma \ref{temam book}. 
%We assume that $\vec{u}_{0} \in {L}^{2}(\mathbb{T}_3)^3)$ where $\vec{u}_{0}$ is divergence-free,  $\vec{f} \in L^{2}([0,T],L^{{{2}{}}}(\mathbb{T}_3)^3)$ and $ \vec{u}$ is a weak solution to (\ref{alpha ns}).
Let $( \alpha_i, \beta_i) $, $i \in I$, be the connected components of $ \mathcal{O}$. Then
\begin{eqnarray}
  \label{leray reformulation}
\displaystyle  \sum_{i \in I}( \beta_i - \alpha_i)^{\frac{1-4\theta}{2}}<\infty
\end{eqnarray}
\end{lemme}
\textbf{Proof.}
Let $( \alpha_i, \beta_i) $ be one of these connected components and let $t \in ( \alpha_i , \beta_i) \subseteq \mathcal{O}$.   Since
$\vec{u}_{} \in  C_{weak}([0,T];\vec{H}) \cap L^{2}([0,T];\vec{V})$, $\vec{u}_{}(t)$ is well defined for every $t \in ( \alpha_i , \beta_i)$ and $t$ can be chosen   such that $\vec{u}_{}(t) \in \vec{V} $. 
 According to Theorem \ref{2deuxieme}, inequality (\ref{maximal time}) and the fact that $ \| \vec{u}(\beta_i)\|_{\vec{V}}=+\infty$,  we have for $t \in ( \alpha_i , \beta_i) $
  \begin{eqnarray*}
  \label{maximal time2}
  \displaystyle \beta_i - t \ge   \frac{1}{  C(\alpha,\theta, \nu, f )}\frac{1}{(1+\|\vec{u}_{}(t)\|_{\vec{V}}^2)^{\gamma -1} },
  \end{eqnarray*}
   where we have used that $\gamma= \frac{3+4\theta}{1+4 \theta}   >1 $. 
  Thus
  \begin{eqnarray*}
  \label{maximal time3}
  \displaystyle \frac{ C(\alpha,\theta, \nu, f )}{(\beta_i - t)^{\frac{1}{\gamma -1}} }  \le  1+\|\vec{u}_{}(t)\|_{\vec{V}}^2.
  \end{eqnarray*}
  Then we integrate on $( \alpha_i , \beta_i) $ to obtain
   \begin{eqnarray*}
  \label{maximal time4}
  \displaystyle { C(\alpha, \theta,\nu, f )}{(\beta_i - \alpha_i)^{\frac{-1}{\gamma -1}+1} }  \le (\beta_i - \alpha_i) +\displaystyle  \int_{\alpha_i}^{\beta_i} \|\vec{u}_{}(t)\|_{\vec{V}}^2 dt,
  \end{eqnarray*}
  Adding all these relations for $ i \in I $ we obtain
   \begin{eqnarray*}
  \label{maximal time5}
  \displaystyle { C(\alpha, \theta, \nu, f )} \sum_{i \in I}{(\beta_i - \alpha_i)^{\frac{-1}{\gamma -1}+1} }  \le T +\displaystyle  \int_{0}^{T} \|\vec{u}_{}(t)\|_{\vec{V}}^2 dt.
  \end{eqnarray*}
 $\hfill\blacksquare$\\

\textbf{Proof of Theorem \ref{main}.}
We set $ \mathcal{S}=\mathcal{S}_\theta(\vec{u})=[0,T]\setminus \mathcal{O}$. We have to prove that the $\frac{1-4\theta}{2}$-dimensional Hausdorff measure of $\mathcal{S} $ is zero. 
Since the Lebesgue measure of $ \mathcal{O}$ is finite ,i.e.
 \begin{eqnarray}
 \label{finitemeasure}
 \displaystyle \sum_{i \in I}(\beta_i -\alpha_i) < \infty,
  \end{eqnarray}
it follows from Lemma \ref{temam book2} that for every $\epsilon >0$  there exist a finite part $I_{\epsilon} \subset I$ such that
 \begin{eqnarray}
\displaystyle  \sum_{i \in I\setminus I_{\epsilon}}(\beta_i -\alpha_i) \le \epsilon
 \end{eqnarray}
 and
   \begin{eqnarray}
   \label{keyone}
\displaystyle  \sum_{i \in I\setminus I_{\epsilon}}(\beta_i -\alpha_i)^{\frac{1-4\theta}{2}} \le \epsilon
 \end{eqnarray}
 Note that $ \mathcal{S} \subset [0,T] \setminus \displaystyle  \bigcup_{i \in I_{\epsilon}}(\alpha_i,\beta_i)$ and the set $  [0,T] \setminus\displaystyle  \bigcup_{i \in I_{\epsilon}}(\alpha_i,\beta_i)$ 
 is the union of finite number of mutually disjoint closed intervals, say $B_j$, for $j=1,...,N$. Our aim now is to  show that the $\diameter B_j \le \epsilon$.
 Since the intervals $(\alpha_i, \beta_i)$ are mutually disjoint, each interval $\displaystyle  (\alpha_i, \beta_i)$, $\displaystyle  i \in I\setminus I_{\epsilon}$, 
 is included in one, and only one, interval $B_j$. We denote by $ I_j$ the set of indice $i$ such that $(\alpha_i, \beta_i ) \subset B_j$. It is clear that 
 $ \displaystyle I_{\epsilon},I_1,...,I_N $ is a partition of $I$ and we have $B_j=(\bigcup_{i \in I_{j}}(\alpha_i,\beta_i))\cup(B_j\cap\mathcal{S})$ for all 
 $j=1,...,N$. It follows from (\ref{finitemeasure}) that
 \begin{eqnarray}
 \label{dimetereepsilon}
 \diameter B_j = \sum_{i \in I_{j}}(\beta_i -\alpha_i) \le \epsilon.
  \end{eqnarray}
 Finally in virtue of  the definition \ref{definitionhaussdrof} and estimates (\ref{dimetereepsilon}), (\ref{keyone})  and since $\displaystyle l^\delta \hookrightarrow l^1$ for all $ 0< \delta < 1$ we have
  \begin{equation}
  \begin{array}{ccccccc}
\displaystyle \mu_{\frac{1-4\theta}{2},\epsilon}(\mathcal{S}) &\le& \displaystyle \sum_{j=1}^{N}(\diameter B_j)^{\frac{1-4\theta}{2}}\\
\displaystyle &\le& \displaystyle \sum_{j=1}^{N}\left(\sum_{i \in I_{j}}(\beta_i -\alpha_i) \right)^{\frac{1-4\theta}{2}}\\
\displaystyle &\le& \displaystyle\sum_{j=1}^{N}\sum_{i \in I_{j}}\left(\beta_i -\alpha_i\right)^{\frac{1-4\theta}{2}}\\
\displaystyle &=& \displaystyle \sum_{i \in I\setminus I_{\epsilon}}\left(\beta_i -\alpha_i \right)^{\frac{1-4\theta}{2}} \le \epsilon.
\end{array}
 \end{equation}
 Letting $\epsilon \rightarrow 0,$ we find $\displaystyle \mu_{\frac{1-4\theta}{2}}(\mathcal{S})=0$ and this  completes the proof. $\hfill\blacksquare$

\end{document}